\documentclass{amsart}
\usepackage{amsthm}
\usepackage{newlfont}
\usepackage{amsmath}
\usepackage{amssymb}
\usepackage{amsfonts}
\usepackage[arrow,curve,matrix,tips]{xy}

\usepackage{hyperref}

\hypersetup{
 pdfauthor = {Olaf~M.~Schn{\"u}rer},
 pdftitle = {Regular Connections on Principal Fiber Bundles over
  the Infinitesimal Punctured Disc},
 pdfsubject = {Mathematics},
 pdfkeywords = {Regular Connection, Gauge Group, Classification},
}

\newtheorem{theorem}{Theorem}[section]
\newtheorem*{theoremone}{{Theorem~\ref{T:standard}}}
\newtheorem*{theoremtwo}{{Theorem~\ref{T:standardnull}}}
\newtheorem{lemma}[theorem]{Lemma}
\newtheorem{corollary}[theorem]{Corollary}
\newtheorem{proposition}[theorem]{Proposition}

\theoremstyle{definition}
\newtheorem{definition}[theorem]{Definition}
\newtheorem{example}[theorem]{Example}
\newtheorem{examples}[theorem]{Examples}

\theoremstyle{remark}
\newtheorem{remark}[theorem]{Remark}
\newtheorem{question}[theorem]{Question}

\newcommand{\dashkomma}{{i.\,e., }}
\newcommand{\dash}{{i.\,e.\ }}

\newcommand{\op}{\operatorname}

\newcommand{\ra}{\rightarrow}
\newcommand{\sra}{\twoheadrightarrow}
\newcommand{\hra}{\hookrightarrow}
\newcommand{\sira}{\stackrel{\sim}{\rightarrow}}

\newcommand{\veps}{\varepsilon}

\newcommand{\Bl}[1]{{\mathbb{#1}}}
\newcommand{\DZ}{\Bl{Z}}
\newcommand{\DN}{\Bl{N}}
\newcommand{\DNplus}{{\Bl{N}^+}}
\newcommand{\DQ}{\Bl{Q}}
\newcommand{\DC}{\Bl{C}}
\newcommand{\Lie}{\op{Lie}}

\newcommand{\LieD}{\mathfrak{d}}

\newcommand{\LieG}{\mathfrak{g}}
\newcommand{\LieB}{\mathfrak{b}}
\newcommand{\LieS}{\mathfrak{s}}
\newcommand{\LieT}{\mathfrak{t}}

\newcommand{\LieU}{\mathfrak{u}}
\newcommand{\LieZ}{\mathfrak{z}}
\newcommand{\LieGL}{\mathfrak{gl}}
\newcommand{\LieSL}{\mathfrak{sl}}

\newcommand{\Z}{\op{Z}}
\newcommand{\Ho}{\op{H}}
\newcommand{\Hom}{\op{Hom}}
\newcommand{\End}{\op{End}}

\newcommand{\Mat}{\op{Mat}}
\newcommand{\Eig}{\op{Eig}}

\newcommand{\id}{\op{id}}
\newcommand{\Eins}{\mathbf{1}}
\newcommand{\Einheit}{\op{E}}
\newcommand{\Jordan}{\op{J}}

\newcommand{\Sym}{\op{Sym}}

\newcommand{\LieGNull}{\LieG^{\op{zero}}}
\newcommand{\LieGLnNull}{\LieGL_n^{\op{zero}}}
\newcommand{\LieSLnNull}{\LieSL_n^{\op{zero}}}
\newcommand{\LieSLzweiNull}{\LieSL_2^{\op{zero}}}

\newcommand{\ol}[1]{{\overline{#1}}}
\newcommand{\ul}[1]{{\underline{#1}}}

\newcommand{\incl}[1]{{{#1}^\ast}}
\newcommand{\mult}[1]{{{(#1\cdot)}^\ast}}

\newcommand{\GL}{\op{GL}}
\newcommand{\SL}{\op{SL}}
\newcommand{\Sp}{\op{Sp}}

\newcommand{\Ker}{\op{ker}}

\newcommand{\Bild}{\op{im}}

\newcommand{\Ad}{\op{Ad}}
\newcommand{\ad}{\op{ad}}

\newcommand{\cek}{\vee}

\newcommand{\dl}{{\,\op{d}}}

\newcommand{\dz}{{\,\op{d}\!z}}

\newcommand{\dlogz}{{\,\op{d}\!\log z}}

\newcommand{\dlogzeta}{{\,\op{d}\!\log\zeta}}

\newcommand{\delz}{{\,\partial_z}}

\newcommand{\Gal}{\op{Gal}}

\newcommand{\inv}{^{-1}}
\newcommand{\can}{\op{can}}

\newcommand{\calO}{\cal{O}}
\newcommand{\calmax}{\mathfrak{m}}

\newcommand{\calG}{\cal{G}}

\newcommand{\calH}{\cal{H}}

\newcommand{\CGL}{\cal{GL}}
\newcommand{\CGLn}{\cal{GL}_n}

\newcommand{\Hganzko}{R^{\DZ \cek}_{H}}

\newcommand{\Rel}{\op{Rel}}

\newcommand{\semidirect}{\rtimes}

\newcommand{\blockdiag}{\op{blockdiag}\;}

\numberwithin{equation}{section}

\begin{document}
\title[Regular Connections on Principal Fiber Bundles]
{Regular Connections on Principal Fiber Bundles over
  the Infinitesimal Punctured Disc} 
\author{Olaf~M.~Schn{\"u}rer}
\address{Mathematisches Institut, Universit{\"a}t Freiburg, 
  Eckerstra{\ss}e\ 1, D-79104 Freiburg, Germany}
\email{schnuerer@pcpool00.mathematik.uni-freiburg.de}

%\date{October 2006, research from 2002/03.}


\begin{abstract}
This paper concerns regular connections on trivial algebraic
$G$-principal fiber bundles over the infinitesimal punctured disc,
where $G$ is a connected reductive linear algebraic group over an
algebraically closed field of characteristic zero.
We show that the pull-back of every regular connection to an
appropriate covering of the infinitesimal punctured disc is gauge
equivalent to a connection of the form
$X z^{-1}\operatorname{d}\! z$ 
for some $X$ in the Lie algebra of $G$.
We may even arrange that the only
rational eigenvalue of $\ad X$ is zero.
Our results allow a classification of regular $\SL_n$-connections up to
gauge equivalence.
\end{abstract}
\maketitle
\tableofcontents
\section{Introduction}\label{S:introduction}
Let $G$ be a linear algebraic group over an algebraically closed field
$k$ of characteristic zero, and let $\LieG$ be its Lie algebra.
The loop group $G((z))=G(k((z)))$ acts as a gauge group  
on the set
$\calG = \LieG \otimes_k k((z)) \dz$ of connections on the trivial
algebraic 
$G$-principal fiber bundle over the infinitesimal punctured disc.
For $G=\GL_n$, this action is given by 
$$g\left[A  \dz\right] 
=\left(g A g\inv + \tfrac{\partial}{\partial z}(g)g\inv\right)
\dz,$$
where $g \in \GL_n((z))$, 
$A \in \LieGL_n \otimes_k k((z))$, 
and
$\tfrac{\partial}{\partial z}$ acts on each entry of the matrix $g$.
If \mbox{$G \subset \GL_n$} is a closed subgroup, this action
induces an action of $G((z))$ on $\calG$.
According to \cite[\S8.2, Definition]{BV}, a connection is regular if it is 
gauge equivalent 
to an element of \mbox{$\LieG \otimes_k k[[z]] z \inv \dz$}.
For a positive integer $m$, we define the  
inclusion $\incl{m}: k((z)) \hra k((z))$ by
$\incl{m}(f(z)) = f(z^m)$. Geometrically, it corresponds 
to an $m$-fold covering of the infinitesimal punctured disc.
We can pull back every connection $A$ to a connection 
$\incl{m}(A)$. 
Now the main results of this article can be formulated.
\begin{theoremone}
Let $G$ be a connected reductive linear algebraic group. There exists
a positive integer $m$, such that for every regular 
connection $A$ the pull-back connection
$\incl{m}(A)$ is gauge equivalent to
$X z\inv \dz$ for a suitable  $X \in \LieG$.
\end{theoremone}
In \cite[\S8.4~(c)]{BV}, a similar statement is given 
for any affine algebraic group over the complex numbers. Its proof,
however, uses analytic methods. 
Here, we give a purely algebraic proof for a connected reductive group.

\begin{theoremtwo}
Let $G$ be a connected reductive linear algebraic group.
For every regular connection $A$, there exists a
positive integer $m$ and an element $X \in \LieG$,
such that the only rational eigenvalue of $\ad X$ is zero and 
the pull-back 
connection $\incl{m}(A)$
is gauge equivalent to 
$X z\inv \dz$.
\end{theoremtwo}

The proofs of Theorems~\ref{T:standard} and \ref{T:standardnull} are
mainly based on the structure theory 
of the group $G$ and its Lie algebra $\LieG$.
We use these theorems and Galois cohomology
in order to get a classification of regular
$\SL_n$-connections up to $\SL_n((z))$-equi\-va\-lence, see
Theorem~\ref{S:slnklassifik} and Remarks~\ref{rem:slnklassifik}, \ref{rem:nice-class}. Our
classification strategy is motivated by \cite{BV}.
Some steps of this strategy can also be applied to other connected
reductive linear algebraic groups, see \cite{OSdiplom}.

This paper is organized as follows.
Using results from \cite{DG}, we define in Section~\ref{S:conngauge}
the action of the gauge group on the space of connections
intrinsically, \dashkomma without choosing a closed embedding of $G$ into
some $\GL_n$.
We repeat the definitions of regular and aligned connections 
given in \cite{BV} and recall that every
regular connection is gauge equivalent to an aligned connection. 
In Section~\ref{S:pullback}, we explain how to pull back connections.
We call two connections related, if they become gauge equivalent in
some covering. Using Steinberg's theorem (cf.~\cite{Steinberg}), we show that for a
connected group $G$ the relatives of a connection $A$ up to gauge
equivalence correspond bijectively
to a set $\Ho^1(K; A)$ defined via Galois cohomology. 
We prove our main Theorems~\ref{T:standard} and \ref{T:standardnull} in Section
\ref{S:standard}. Section~\ref{S:dmoduln} contains the classification
of regular $\GL_n$-connections up to $\GL_n((z))$-equivalence. This is
classical, see for example \cite[\S3]{BV}. 
We explain the relation
between $n$-dimensional (Fuchsian) $D$-modules (see \cite{Manin}) and (regular)
$\GL_n$-connections. If we translate the classification into the
language of $D$-modules, we obtain a result of Manin (\cite{Manin}).
The results of the previous sections are used in Section
\ref{S:klassifikation} in order to classify regular
$\SL_n$-connections up to relationship and up to gauge equivalence.
These classifications use an explicit description of the semisimple
conjugacy classes in the centralizer $\Z_{\GL_n}(X)$, where $X$ is an
element of the Lie algebra $\LieGL_n$. We include this description in
Appendix \ref{App:cent}.
The definition of a Fuchsian connection in Section~\ref{S:fuchszshg}
sounds more natural than the definition of a regular connection. 
It is based upon the notion of a Fuchsian $D$-module given in \cite{Manin}.
We show that regular connections are Fuchsian, and that for
$G=\GL_n$ and $G=\SL_n$, Fuchsian connections are regular.

This paper is a condensed version of my diploma thesis
\cite{OSdiplom}, written in Freiburg in 2002/2003. 
I would like to thank Wolfgang Soergel.
He taught me a lot of the mathematics I know.

\section{Connections and Gauge Group}\label{S:conngauge}
\subsection{Conventions}\label{SS:conventions}
We fix an algebraically closed field $k$ of characteristic zero
and write $\otimes=\otimes_k$, $\Hom=\Hom_k$, and so on.
If we discuss vector spaces, Lie algebras, linear
algebraic groups, or the like, we always mean the corresponding
structures defined over $k$.
We denote the positive integers by $\DNplus$ and define $\DN = \DNplus
\cup \{0\}$.

Let $\calO=k[[z]]$ be the ring of formal power series over $k$, with
maximal ideal $\calmax = zk[[z]]$ and the induced
$\calmax$-topology. The quotient field of 
$\calO$ is the field $K=k((z))$ of Laurent series over $k$.
Let $D= k((z)) \delz$ denote the continuous $k$-linear derivations
from $K$ to $K$, where we abbreviate $\delz = \tfrac{\partial}{\partial z}$.
Define $\Omega=\Hom_K(D,K) = k((z))\dz$, where $\dz$ is dual to
$\delz$, \dashkomma $\dz(\delz)=1$.

\subsection{Action of the Gauge Group on the
  Space of Connections}\label{SS:opergauge}

We collect some results from \cite[II,~\S4]{DG}, 
in particular a definition of the Lie algebra.
These will enable us to define the action of the
gauge group on the space of connections in an intrinsic way.
The reader who is not interested in this intrinsic definition may
use Equation~\eqref{Eq:actionGLn} in
Example~\ref{Bsp:GLnOp} as the definition for closed subgroups of
$\GL_n$. He may check that this is well defined and does not depend on
the closed embedding. 

Let $G$ be a linear algebraic group. We consider $G$ as an affine
algebraic group scheme.
Suppose that $R$ is a $k$-algebra. We denote the algebra of dual
numbers by $R[\veps]=R\oplus R\veps$. 
By applying the group functor $G$ to 
the unique $R$-algebra homomorphism from $R$ to $R[\veps]$, we
consider $G(R)$ as a subgroup of $G(R[\veps])$.
Define the $R$-algebra homomorphism
$p:R[\veps]\ra R$ by $p(\veps) = 0$. 
As explained in \cite{DG}, the kernel $\LieG(R)$ of the group homomorphism
$G(p):G(R[\veps]) \ra G(R)$ is endowed with a structure of
Lie algebra over $R$. 
The Lie algebra $\LieG = \LieG(k)$ is canonically isomorphic to the
standard Lie algebra of the linear algebraic group $G$.
The obvious inclusion $k \hra R$ induces a homomorphism of Lie algebras
$\LieG \ra \LieG(R)$. Tensoring with $R$ gives a
canonical homomorphism 
$R\otimes\LieG \ra \LieG(R)$
of Lie algebras over $R$. This is an isomorphism, 
as our group scheme $G$ is locally algebraic. Hence we identify
$R\otimes\LieG = \LieG(R)$.  

We define an action of the
gauge group $G(K)$ on the set 
$$  \calG 
= \Hom_K(D, \LieG(K))$$
of connections
on the trivial algebraic
$G$-principal fiber bundle over the infinitesimal punctured disc.
Every derivation $s\in D$ gives rise to a $k$-algebra homomorphism
$\hat{s}: K \ra K[\veps]$ defined by $\hat{s}(f) = f + s(f)\veps$.
We get a group homomorphism
$ \hat{s}=G(\hat{s}) : G(K) \ra G(K[\veps])$.
For $g\in G(K)$, define $\dot{g}:  D  \ra  G\left(K[\veps]\right)$ by
$\dot{g}(s)= \hat{s}(g)$.
For all $g\in G(K)$ and $s\in D$, we deduce from
$p\circ \hat{s}=\id_K$ that 
$\dot{g}(s)g\inv$
is an element of the kernel of $G(p)$, \dashkomma an element of
$\LieG(K)$.

\begin{proposition}\label{P:gaugeaction}
The map 
$$G(K) \times \calG  \ra  \calG,\quad
(g,A)  \mapsto g[A] = (\Ad g) \circ A + (\cdot {g\inv})\circ \dot{g},$$
defines an action of the gauge group $G(K)$ on the space of connections
$\calG$. Here, $(\cdot h)$ denotes right multiplication by a group
element $h$.

If $f:G \ra H$ is a homomorphism of linear algebraic groups, we have
$f_{\ast}(g[A]) = f(g)[f_{\ast}(A)]$, where $f_{\ast}:\calG \ra \calH$
is the obvious map.
\end{proposition}

\begin{proof}
We have to show that, for $g \in G(K)$ and $A\in\calG$, the map
$g[A]$ from $D$ to $\LieG(K)$ is $K$-linear.
In order to do this, use the definition in \cite{DG} of the $K$-vector
space structure on $\LieG(K)$. The explicit arguments can be found
in \cite{OSdiplom}.
The remaining claims are obvious.

There is another way of proving the $K$-linearity of the map $g[A]$:
By choosing a closed embedding $G \hra \GL_n$, we reduce to the case
$G=\GL_n$.
Then an easy calculation based on 
Equation~\eqref{Eq:actionGLn} in
Example~\ref{Bsp:GLnOp} shows the $K$-linearity.  
\end{proof}

\begin{definition}
Let $H\subset G(K)$ be a subgroup. Two connections $A$, $B \in \calG$
are {\bf $H$-equivalent}, if there is an element $h \in H$
such that $h[A] = B$. They are {\bf gauge equivalent}, if they
are $G(K)$-equivalent.
\end{definition}
Note that canonically
\begin{equation*}
  \calG=\Hom_K(D, \LieG(K)) = \LieG(K)\otimes_K \Hom_K(D, K) =
  \LieG(K) \otimes_K \Omega. 
  % = \LieG \otimes \Omega. 
\end{equation*}
Since $\Omega =K \dz$, each
element of $\LieG(K) \otimes_K \Omega$ can be written uniquely in the
form $A \otimes \dz$ with $A \in \LieG(K)$. We abbreviate $A \dz = A
\otimes \dz$ and write $\calG = \LieG(K)\dz$ accordingly.
We define $\dlogz = z\inv \dz \in \Omega$ and write similarly $A \dlogz$
and $\calG = \LieG(K)\dlogz$.
If $B \in \calG = \Hom_K(D, \LieG(K))$ is a connection, 
we have $B = B(\delz)\dz = B(z\delz) \dlogz$.

We use the abbreviations 
$G((z))$ (resp.\ $\LieG((z))$, $\LieG[z]$) for $G\left(k((z))\right)$
(resp.\ $\LieG\left(k((z))\right)$, $\LieG\left(k[z]\right)$), and so on.

For a connection $A \dlogz \in \LieG((z))\dlogz$ and 
a gauge transformation $g \in G((z))$, we get
\begin{equation}\label{Eq:actionG}
g[A \dlogz] = \left((\Ad g)(A) + \widehat{z\delz}(g)g\inv\right)\dlogz.
\end{equation}

\begin{example}\label{Bsp:GLnOp}
Consider the case $G=\GL_n$.
Let 
$A \in \Mat_n\left(k((z))\right)=\LieGL_n((z))$
and 
\mbox{$g \in \GL_n((z))$} be given. 
As $\widehat{z\delz}(g)=g + z\delz (g)\veps$ is satisfied in
$\GL_n(k((z))[\veps])$,  
we get 
\begin{align}\label{Eq:actionGLn}
g\left[A \dlogz\right] 
=& \left(g  A  g\inv + z\delz(g)g\inv\right) \dlogz.
\end{align}
\end{example}

\subsection{Alignement of Regular Connections}\label{SS:alignregconn}
Let $X\in\LieG$ and $\lambda \in k$. We denote the eigenspace of $\ad
X$ corresponding to $\lambda$ by
$\Eig(\ad X; \lambda)=\ker(\ad X - \lambda)$. Let $X=
X_{{\op{s}}}+X_{\op{n}}$ be the Jordan decomposition in $\LieG = \Lie G$.

\begin{definition}(cf.~\cite[\S8.2, Definition, and \S8.5, Definition]{BV})
The elements of $\LieG[[z]]\dlogz$ are called 
{\bf connections of the first kind}.
A connection is {\bf regular}, if it is
gauge equivalent to a connection of the first kind.
A connection of the first kind 
$A = \sum_{r \in \DN}A_rz^r \dlogz$ 
is {\bf aligned}, if 
$A_r \in \Eig (\ad {A_{0,{\op{s}}}};r)$ 
for all $r \in \DN$.
\end{definition} 
\begin{theorem}[{\cite[\S8.5, Proposition]{BV}}]\label{T:regaus}
Every regular connection is gauge equivalent to an
aligned connection.
\end{theorem}
For a worked out proof, that uses the exponential map as defined in
\cite[II,~\S6,~3]{DG}, see \cite{OSdiplom}.

\section{Relatives and Galois
  Cohomology}\label{S:pullback}
\subsection{Pull-Back of Gauge Transformations and Connections}\label{SS:pullbackdiffformsandconn}
\begin{definition}
  Let $E$ and $F$ be fields, $\phi:E\ra F$ a map, $V$
  an $E$-vector space, and $W$ an $F$-vector space.
By a $\phi$-{\bf linear} map $f:V\ra W$ we mean a group homomorphism
$f:V\ra W$ such that $f(ev) = \phi(e)f(v)$ for all $e \in E$, $v \in V$.
\end{definition}

\begin{definition}
For $m \in \DNplus$ we define 
the field extension $\incl{m}: K \hra K=M$ by
$\incl{m}(f(z)) = f(z^m)$.
\end{definition}

Let $G$ be a linear algebraic group, let $l$, $m \in \DNplus$, and let
$\phi:\incl{l} \ra \incl{m}$ be a morphism of field extensions, 
\dash\ $\phi:K\ra K$ is a ring homomorphism satisfying \mbox{$\phi \circ
  \incl{l} = \incl{m}$}.
The map $\phi=G(\phi): G(K) \hra G(K)$ is
called $\phi$-pull-back of gauge transformations.
There is a unique $\phi$-linear map
\mbox{$\Phi: \Omega \ra \Omega$} 
such that $\dl \circ \phi = \Phi \circ \dl$,
where the derivation $\dl:  K  \ra  \Omega$ is defined by
$\dl x (\delta) = \delta (x)$, for $x \in K$, $\delta \in D$. 
We denote this injective map by $\phi$ and call it $\phi$-pull
back of $1$-forms. 
This map $\Phi$ and the $\phi$-linear homomorphism of Lie algebras
$\LieG(\phi):\LieG(K) \hra \LieG(K)$
induce an injective $\phi$-linear map
$$\LieG(\phi) \otimes \Phi:\calG=\LieG(K)\otimes_K \Omega \hra
\calG=\LieG(K)\otimes_K \Omega.$$
We denote this map simply by $\phi$ and call it $\phi$-pull-back of
connections. 

\begin{example}
In particular, for $l=1$ and $\phi= \incl{m}$, we have
$\incl{m}(\dlogz)=m\dlogz$,
and therefore, for $A(z) \in \LieG((z))$,
\begin{equation}\label{Eq:pullbacki}
\incl{m}(A(z) \dlogz)= m A(z^m)\dlogz.
\end{equation}
\end{example}

\begin{proposition}\label{P:phiswitch}
Under the above assumptions we have
$\phi(g)[\phi(A)]=\phi(g[A])$ for all $g \in G(K)$, $A \in \calG$.
\end{proposition}
\begin{proof}
  Use a closed embedding $G \subset \GL_n$ and
  Equation~\eqref{Eq:actionGLn} in Example~\ref{Bsp:GLnOp}.
  For a more intrinsic proof we refer to \cite{OSdiplom}.
\end{proof}

\subsection{Connections and Galois Cohomology}\label{SS:formsandgalois}
For $m \in \DNplus$ let $\Gamma_m= \Gal(\incl{m})$ be the Galois group
of the field extension $\incl{m}:K \hra K=M$.
The map
$\Gamma_m \times \calG \ra \calG$, 
$(\sigma, A) \mapsto \sigma(A)$,
given by the pull-back of connections,
is an 
action of $\Gamma_m$ on $\calG$.
\begin{lemma}\label{L:invariantzshg}
  Let $m \in \DNplus$.
  A connection $C$ is 
  \mbox{$\Gamma_m$-invariant} if and only if there is a 
  connection $D$ such that $C = \incl{m}(D)$.
\end{lemma}
\begin{proof}
  From $\sigma \circ \incl{m} = \incl{m}$, we deduce that
  \mbox{$\sigma(\incl{m}(D))=\incl{m}(D)$} for all $\sigma \in \Gamma_m$ and all
  connections $D$.
  Conversely, suppose that $C(z)\dlogz \in
  \calG$ is $\Gamma_m$-invariant. 
  For $\omega$ a primitive
  $m$-th root of unity,
  we define $\mult{\omega} \in \Gamma_m$ by
  $\mult{\omega}(f(z))= f(\omega z)$.
  From $\mult{\omega}(C(z)\dlogz) = C(\omega z) \dlogz$ we see that
  $C(z) = B(z^m)$ for some $B \in \LieG((z))$. 
  Equation~\eqref{Eq:pullbacki} implies that the connection 
  $D=m\inv B (z)\dlogz$ satisfies $\incl{m}(D) = C (z)\dlogz$.
\end{proof}

\begin{definition}
  Let $C$ be a connection, $m \in \DNplus$. An $m$-{\bf form} of $C$ is a
  connection $D$ such that $ \incl{m}(D)$ and $C$ are gauge equivalent.
\end{definition}

Fix a $\Gamma_m$-invariant connection $A$.
The action of the Galois group $\Gamma_m$ on $K=M$ defines an
action on $G(M)$ by group automorphisms. 
By Proposition~\ref{P:phiswitch}, this action restricts
to an action on
the stabilizer 
$G(M)_A$ of the $\Gamma_m$-invariant connection $A$.
We define a map
$$p=p(A):\left\{\text{$m$-forms of $A$} \right\} \ra \Ho^1(\Gamma_m;
G(M)_A), \quad B \mapsto p^B,$$
as follows. Given an $m$-form $B$ of $A$, choose $b \in G(M)$
such that $b[\incl{m}(B)]=A$.
Then the map $p^b:\Gamma_m \ra G(M)_A$ defined by 
$p_\sigma^b = b\sigma(b\inv)$ is a $1$-cocycle, 
and its cohomology class
$p^B= [p^b]$ does not depend on the choice of $b$.

\begin{theorem}\label{T:H1Kform}
  For $m \in \DNplus$, consider the field extension $\incl{m}:K\hra
  K=M$. If $A$ is a $\Gamma_m$-invariant
  connection, the map $p=p(A)$ defined above descends
  to an injection
  $$\ol{p}=\ol{p(A)}:\left\{\text{$m$-forms of $A$} \right\}/G(K) \hra
  \Ho^1(\Gamma_m; G(M)_A).$$
  If $\Ho^1(\Gamma_m; G(M)) = \{1\}$ 
  this map $\ol{p}$ is bijective. 
\end{theorem}
\begin{proof}
  Let $\Gamma = \Gamma_m$.
  Let $B$ and $C$ be gauge equivalent $m$-forms of $A$. Let
  $g \in G(K)$ with $B=g[C]$, and let $b \in G(M)$ with $b[\incl{m}(B)] =
  A$. Consequently, we have $(b\incl{m}(g))[\incl{m}(C)]= A$. Now
  $$p_{\sigma}^{b\incl{m}(g)}=b\incl{m}(g)\sigma(\incl{m}(g)\inv
  b\inv) = b \sigma(b\inv) = p_\sigma^b$$
  implies that $p^B=p^C$. Thus our map $p$ descends.

  Let $B$ and $C$ be $m$-forms of $A$. Choose $b$, $c \in G(M)$ with
  $b[\incl{m}(B)]= c[\incl{m}(C)]= A$. If $p^b$ and $p^c$ become equal
  in $\Ho^1(\Gamma, G(M)_A)$, there is an element $f\in G(M)_A$ such
  that $b\sigma(b\inv)= fc\sigma(c\inv f\inv)$ for all $\sigma \in
  \Gamma$. We deduce that the element $c\inv f\inv b$ is
  $\Gamma$-invariant and hence equal to $\incl{m}(h)$ for some $h \in
  G(K)$. But then
  \begin{equation*}
    \incl{m}(h[B]) = \incl{m}(h)[\incl{m}(B)] = c\inv f\inv b
    [\incl{m}(B)] 
    = \incl{m}(C)
  \end{equation*}
  and $h[B]=C$ since $\incl{m}$ is injective. So $\ol{p}$ is injective.
 
  Suppose that $\Ho^1(\Gamma; G(M)) = \{1\}$. We now prove that $p$
  is surjective.
  Let $a\in \Z^1(\Gamma; G(M)_A)$ be a $1$-cocycle. By assumption, $a$
  considered as an element of 
  $\Z^1(\Gamma; G(M))$
  is cohomologous to the trivial $1$-cocycle, \dashkomma there is $g \in
  G(M)$ such that 
  $a_\sigma = g \sigma(g\inv)$ for all $\sigma \in \Gamma$.
  Let $C = g\inv [A]$. For any $\sigma \in \Gamma$, we get
  $$\sigma(C) =\sigma(g\inv)[A] = (g\inv a_\sigma)[A] = C.$$
  According to Lemma~\ref{L:invariantzshg}, there is a
  connection $D$ 
  such that $\incl{m}(D) =C$. 
  As $g[\incl{m}(D)]=A$, we have $p^D=[p^g] = [a]$.
\end{proof}

\begin{remark}\label{R:H1triv}
  In the cases  $G=\GL_n$, $G=\SL_n$ or $G=\Sp_{2n}$
  we know that $\Ho^1(\Gamma_m; G(M)) = \{1\}$, according to \cite[X,~\S1
  and \S2]{Serre}. 

  We now explain that, more generally, for every connected linear
  algebraic group $G$ we have $\Ho^1(\Gamma_m; G(M)) = \{1\}$.
  It is well known \cite[IV,~\S2, Proposition~8 and Corollary]{Serre} that
  the union 
  $\ol{K}=\bigcup_{m\in\DNplus}k((z^{1/m}))$
  is an algebraic closure of $K=k((z))$.
  Thus, every algebraic extension of $K$ ramifies,
  \dashkomma $K$ itself is the maximal 
  unramified extension of $K$. 
  According to \cite[Theorem~12]{Lang} and \cite[II,~\S3.2
  Corollary]{Serregaloiscoho}, 
  the maximal unramified extension of a field 
  that is complete with respect to a discrete valuation and
  that has perfect residue class field has 
  dimension $\leq 1$. 
  Consequently, we have $\dim K \leq 1$.

  Let $G$ be a connected group.
  Since $\ol{K} \otimes k[G]$ is an integral domain, $G(\ol{K})$
  is a connected linear algebraic group over $\ol{K}$,
  and $K \otimes k[G]$ defines a $K$-structure on
  $G(\ol{K})$ in the sense of \cite{Springerneu}.
  By \cite[Theorem~17.10.2]{Springerneu} 
  (cf.~\cite{Steinberg} and
  \cite[III,~\S2.3, Theorem~1']{Serregaloiscoho}), 
  we know that $$\Ho^1(\Gal(\ol{K}/K); G(\ol{K})) = \{1\}.$$

  For $m \in \DNplus$, we view the field extension $\incl{m}: K \hra
  K=M$ as a subextension of $K \subset \ol{K}$ via the embedding
  $K=M \hra \ol{K}$, $f(z) \mapsto f(z^{1/m})$.
  As the canonical map from 
  $\Ho^1(\Gamma_m; G(M))$
  to the direct limit
  $$\Ho^1(\Gal(\ol{K}/K); G(\ol{K})) = 
    \varinjlim_{m\in \DNplus}\, \Ho^1(\Gamma_m; G(M))$$
  is injective, we deduce that $\Ho^1(\Gamma_m; G(M)) = \{1\}$. 
\end{remark}

An easy calculation, left to the reader, Theorem~\ref{T:H1Kform}, and
Remark~\ref{R:H1triv} yield

\begin{proposition}\label{P:H1fieldinclusion}
  Let $l$, $m \in \DNplus$ with $d=m/l \in \DNplus$. Consider the
  commutative diagram of field extensions
  \begin{equation*}
    \xymatrix@C10pt{
      K=L \ar[rr]^-{\incl{d}} && K=M\\
      & K \ar[lu]^-{\incl{l}} \ar[ru]_-{\incl{m}}
    }
  \end{equation*}
  If $A$ is a connection, 
  the injection
$$ \incl{d}:\Z^1(\Gamma_l; G(L)_{\incl{l}(A)}) \hra
\Z^1(\Gamma_m; G(M)_{\incl{m}(A)}),$$
defined by $(\incl{d}(p): \sigma \mapsto \incl{d}(p_{\sigma|_{L}})$,
induces an injection $\incl{d}=\incl{(m/l)}$ on cohomology.
  Furthermore, the diagram 
  $$
  \xymatrix@C+10pt@R-10pt{
    {\left\{\text{$l$-forms of $\incl{l}(A)$} \right\}/G(K)
      \ar@{^{(}->}[r]^-{\ol{p(\incl{l}(A))}}} \ar@{^{(}->}[d] & \Ho^1(\Gamma_l; G(L)_{\incl{l}(A)})
    \ar@{^{(}->}[d]^-{\incl{(m/l)}} \\
    {\left\{\text{$m$-forms of $\incl{m}(A)$} \right\}/G(K)
      \ar@{^{(}->}[r]^-{\ol{p(\incl{m}(A))}}} & \Ho^1(\Gamma_m; G(M)_{\incl{m}(A)})
  }
  $$
  commutes.
  If $G$ is connected, the horizontal maps in this diagram are
  bijective.
\end{proposition}

\begin{definition}
  Let $A$ be a connection. The punctured sets
  $\Ho^1(\Gamma_m; G(M)_{\incl{m}(A)})$, for $m \in \DNplus$, 
  together with the maps $\incl{(m/l)}$,
  for $l$, $m \in \DNplus$ with $l$ divides $m$,
  form a directed system. 
  We denote its direct limit by $\Ho^1(K; A)$. 
\end{definition}

\begin{definition}\label{d:relatives}
  Two connections $A$ and $B$ are {\bf related} if there is $m \in
  \DNplus$ such that $\incl{m}(A)$ and $\incl{m}(B)$ 
  are gauge equivalent. This defines an equivalence relation on the
  set of connections, and we define $\Rel(A)$ to be the equivalence
  class containing $A$. The elements of $\Rel(A)$ are called the
  {\bf relatives} of $A$. 
\end{definition}

\begin{proposition}\label{P:verwandteh1}
  Let $A$ be a connection. The map $\Rel(A)/G(K) \hra \Ho^1(K; A)$
    that is induced by the maps $\ol{p(\incl{m}(A))}$, for $m \in \DNplus$,
    is injective. It is bijective if $G$ is a connected group.
\end{proposition}

\begin{proof}
  The set of relatives of $A$ is the union, \dashkomma the direct
  limit, of all $m$-forms of all connections
  $\incl{m}(A)$ for $m \in \DNplus$. Now use Proposition
  \ref{P:H1fieldinclusion}. 
\end{proof}

\section{Transforming Regular Connections}\label{S:standard}
\subsection{Transforming Regular Connections to Standard
  Form}\label{SS:standard}
\begin{definition}
Let $G$ be a linear algebraic group. 
Elements of $\LieG \dlogz$ are called 
{\bf connections in standard form}.
\end{definition}

\begin{theorem}\label{T:standard}
Let $G$ be a connected reductive linear algebraic group.
There exists a positive integer $m \in \DNplus$,
such that for every regular connection $A$ the pull-back 
connection $\incl{m}(A)$ is gauge equivalent  
to a connection in standard form.
\end{theorem}
\begin{corollary}\label{c:regular-related-standard}
  If $G$ is connected reductive, each regular connection is related to
  a connection in standard form. 
\end{corollary}
\begin{proof}
  If the pull-back $\incl{m}(A)$ of a connection $A$ is gauge equivalent
  to $X \dlogz$ for some $X \in \LieG$, then Equation
  \eqref{Eq:pullbacki} shows that $A$ is related to $m\inv X \dlogz$.
\end{proof}
\begin{examples}\label{Ex:standard}
  1. For $G=\GL_n$, we can choose $m=1$, as one can see from the choice
  of $m$ in the proof of Theorem~\ref{T:standard}.

  2. Let $k=\DC \subset K = \DC((z))$ and $G=\SL_2$.
    By the proof of Theorem~\ref{T:standard}, we are able to choose $m = 2$.
    But one has to choose $m>1$.
    % a proper field extension $M$ of $K$: 
    Given
    $n \in \DNplus$, we define
    the regular aligned $\SL_2$-connection
    \begin{equation*}
      A_{n}= 
      \begin{bmatrix}
        \tfrac{n}{2}   & z^n \\
        0 & -\tfrac{n}{2}
      \end{bmatrix}
      \dlogz.
    \end{equation*}
    Thanks to \cite[\S8.2 Example]{BV} we know that $A_{n}$ is
    $\SL_2(K)$-equivalent to a connection in
    standard form if and only if $n$ is even.
    If $k$ is algebraically closed of characteristic $\not=2$, the
    same statement can be verified by direct computation.
\end{examples}

We now recall some results from \cite{Springerneu} and
\cite{SpringerAMS}. 
Let $G$ be a connected reductive linear algebraic group,
$T \subset G$ a maximal torus and 
$(\mathcal{X}(T), \mathcal{R}, \mathcal{X}^\cek(T),
\mathcal{R}^\cek)$ the associated root datum.
The derived group $G'=(G,G)$ is semisimple.
Let $T'$ be the subgroup of $T$ generated by the images of all
$\alpha^\cek$, $\alpha \in 
\mathcal{R}$. To the maximal torus $T'$ in $G'$ we associate
the root datum $(X(T'), R, X^\cek(T'), R^\cek)$.
The restriction map $\mathcal{X}(T) \sra X(T')$ induces a canonical
identification $\mathcal{R}=R$.

On the Lie algebra level, 
the reductive Lie algebra $\LieG$ is the direct sum of its center
$\LieZ$ and the semisimple Lie algebra $[\LieG,\LieG]=\LieG'$.
We have $\LieT = \LieZ \oplus \LieT'$.
We associate to the Cartan subalgebra $\LieT'$ in $\LieG'$
the roots $\ul R$
in ${\LieT'}^\ast=\Hom(\LieT',k)$.
Let $\ul{R}^\cek \subset \LieT'$ denote the coroots.
There are canonical identifications $R=\ul{R}$ and
$R^\cek=\ul{R}^\cek$.

By taking the derivative at the unit element, every root $\alpha\in R
=\mathcal{R}$  can be considered as an element of $\LieT^\ast$.
For $H\in\LieT$ and $\alpha\in R$, we have
$\langle\alpha,H\rangle=\langle\alpha,H'\rangle$ 
if we decompose $H$ as $H=H''+H'\in \LieZ \oplus \LieT'=\LieT$.

Suppose that $H\in\LieT$ is arbitrary. Let
$$R_H^\DZ=\{\alpha \in R \mid \langle \alpha, H \rangle
\in \DZ\} \subset R $$
denote the roots integral on $H$, and let
$$\Hganzko=\{\alpha^{\cek} \mid \alpha \in R_H^\DZ\} \subset R^{\cek}$$
denote the corresponding coroots.
Define  
$$V_H=\DQ R_H^\DZ \subset {\LieT'}^\ast
\quad \text{and}\quad
V_H^\cek=\DQ \Hganzko %\subset {\LieT'}_\DQ %=\DQ \,\ul{Q}^\cek 
\subset {\LieT'}.$$ 
The roots integral on $H$ are a root system $R_H^\DZ$ in $V_H$. 
The canonical map 
$$ V_H^\cek \sira  \Hom_\DQ(V_H, \DQ), \quad
  \lambda  \mapsto  \lambda|_{V_H},$$
is an isomorphism of $\DQ$-vector spaces. Therefore, we identify 
$V_H^\cek = \Hom_\DQ(V_H, \DQ)$.

To the dual root system
$\Hganzko$ in $V_H^\cek$, we associate
the root lattice $Q(\Hganzko) = \DZ \Hganzko$.
It is a subgroup of finite index in the
weight lattice 
$$P(\Hganzko) = \{x\in V_H^\cek \mid \langle \alpha, x \rangle \in \DZ
\quad \text{for all $\alpha \in R_H^\DZ$}\}.$$
From $Q(\Hganzko) \subset X^\cek(T') \subset
\mathcal{X}^\cek(T)$,
we obtain that
$$\left| P(\Hganzko) / Q(\Hganzko)\right| \in 
\{ m \in \DNplus \mid m \cdot P (\Hganzko) \subset
\mathcal{X}^\cek(T)\}.$$
In particular, the set on the right hand side is not empty.
We denote its minimum by $m_H$.
Note that the sets $\{\Hganzko \mid H \in \LieT\}$ and 
$\{m_H \mid H \in \LieT\}$ are finite.

\begin{proposition}\label{P:olmexist}
Let  
$m= \op{lcm}\{m_H \mid H \in \LieT\}$
be the least common multiple of all $m_H$.
For every $H\in \LieT$, there is a cocharacter $\psi \in
\mathcal{X}^\cek(T)$, such that
$ \langle \alpha, \psi \rangle = m \langle \alpha, H \rangle$
for all $\alpha \in R_H^\DZ$.
\end{proposition}
This proposition is a consequence of 
\begin{lemma}\label{L:mHtuts}
Let $H \in \LieT$.
There is a cocharacter $\phi \in \mathcal{X}^\cek(T)$, such that
$\langle \alpha, \phi \rangle = m_H \langle \alpha, H \rangle$
for all $\alpha \in R_H^\DZ$.
\end{lemma}

\begin{proof}[Proof of Lemma~\ref{L:mHtuts}]
Let $B_H = \{\alpha_1,\dots,\alpha_l\}$,
where $\alpha_i \neq \alpha_j$ for $i\neq j$,
be a basis of the root system $R_H^\DZ$ in $V_H$.
Let $\{\varpi_1^\cek,\dots, \varpi_l^\cek\}$ be the
basis dual to $B_H$ in $V_H^\cek$, characterized by  
$\langle \alpha_i,\varpi_j^\cek\rangle =
\delta_{ij}$ for all $i$, $j \in \{1,\dots,l\}$. 
Define $\tau = \sum_{i=1}^l \langle \alpha_i, H\rangle \varpi_i^\cek$. 
As $\tau$ is an element of 
$P(\Hganzko)=\bigoplus_{i=1}^l \DZ\varpi_i^\cek$, the definition
of $m_H$ shows that
$\phi = m_H \tau$ is an element of $\mathcal{X}^\cek(T).$
We have $\langle \alpha_j, \phi \rangle = m_H \langle \alpha_j,
H\rangle$
for all $j \in \{1,\dots,l\}$. 
As $B_H$ is a basis of $R_H^\DZ$, the claim follows. 
\end{proof}

\begin{proof}[Proof of Theorem~\ref{T:standard}]
Let $m \in \DNplus$ be the smallest positive integer that satisfies the
following condition.
$$\text{For all $H\in \LieT$, there is 
$\psi \in \mathcal{X}^\cek(T)$, such that 
$\langle \alpha, \psi\rangle = m \langle\alpha, H\rangle$
for all
$\alpha \in R_H^\DZ$.}$$
By Proposition~\ref{P:olmexist}, this is well defined.
Let $A$ be a regular connection. By Theorem~\ref{T:regaus}, we
may assume that $A$ is aligned. 
Therefore, there are
elements $N\in\DN$ and $A_0,\dots,A_N\in\LieG$ such that
\begin{align*}
A        =& \left( A_0+A_1z+\dots+A_Nz^N \right)\dlogz \quad\text{and}\\
A_r    \in& \;\Eig (\ad {A_{0,\op{s}}}; r) \quad\text{for all $r
  \in\{0,\dots,N\}$.}
\end{align*}
Decompose $A_r={A''}_r+A{'}_r \in \LieZ\oplus\LieG'=\LieG$ for $r
\in\{0,\dots,N\}$. 
Since all Cartan subalgebras of a semisimple Lie algebra are conjugate
under the adjoint group, we find 
$x\in G'$ such that
$(\Ad x) ({A'}_{0,{\op{s}}}) \in \LieT'$.
As $x \in G'(k) \subset G(K)$, we have $\widehat{z\delz}(x)=x$, \dashkomma
$\widehat{z\delz}(x)x\inv=0$ in $\LieG(K)$. 
Using Equation \eqref{Eq:actionG}, we see that
$$ x[A] = \big((\Ad x)(A_0)+(\Ad x)(A_1)z+\dots+(\Ad x)(A_N)z^N\big) \dlogz$$
is also aligned. 
Thus, by replacing $A$ by $x[A]$, we may assume that 
${A'}_{0,{\op{s}}}\in\LieT'$.
We define $H = A_{0,{\op{s}}}$ and
have $H={A''}_0+{A'}_{0,{\op{s}}} \in \LieZ \oplus \LieT'=\LieT.$
As $A$ is aligned, we know that
$A_r \in \Eig (\ad H; r)$ for all \mbox{$r\in\{0,\dots,N\}$}.

By the definition of $m$, we find 
$\psi \in \mathcal{X}^\cek(T)$ such that 
$$\langle \alpha, \psi \rangle = m \langle\alpha, H \rangle
\quad\text{for all $\alpha \in R_H^\DZ$.}$$
Note that $\psi:G_{\text{m}} \ra T$ is a homomorphism of group schemes from
the multiplicative group to our torus.
Defining \mbox{$t = \psi(z\inv) \in T(K)$}, we get
$$\alpha(t)=z^{-\langle \alpha, \psi \rangle} 
= z^{-m\langle\alpha, H\rangle}
\quad\text{for all
  $\alpha \in R_H^\DZ$}.$$

We claim that $t[\incl{m}(A)]$ is an element of $\LieG \dlogz$.
Equations~\eqref{Eq:actionG} and \eqref{Eq:pullbacki} yield that
%\begin{align*}
% t[\incl{m}(A)]
%   =&\Big(m(\Ad t) \left(A_0\right) + m(\Ad t) \left(A_1 z^m\right)+\dots
%   + m(\Ad t) \left(A_N z^{mN}\right) \\
%&\quad\quad\quad + \widehat{z\delz}(t)t\inv \Big)
%   \dlogz.
%\end{align*}
\begin{equation*}
 t[\incl{m}(A)]
   =\Big(m(\Ad t) \left(A_0\right) %+ m(\Ad t) \left(A_1 z^m\right)
   +\dots
   + m(\Ad t) \left(A_N z^{mN}\right) + \widehat{z\delz}(t)t\inv \Big)
   \dlogz.
\end{equation*}
Lemma~\ref{L:Adt} implies that we have $(\Ad t) (A_rz^{mr}) = A_r$ 
for all $r \in \{0, \dots, N\}$. 
As $\widehat{z\delz}(t)t\inv$ is an element of $\LieT$,
by Lemma~\ref{L:torustrafo}, 
our claim follows. 
\end{proof}
\begin{lemma}\label{L:Adt}
Let $H \in \LieT$, $m \in \DNplus$, and 
$t\in T(K)$ be such that
$\alpha(t)=z^{-m\langle \alpha, H \rangle}$ 
for all $\alpha \in R_H^\DZ$.
Then, for all $r \in \DN$ and all $B \in \Eig(\ad H; r)$, we have 
$(\Ad t) (Bz^{mr}) = B$.
\end{lemma}

\begin{proof}
We decompose 
$B = B_0+\sum_{\alpha\in R}B_{\alpha}\in 
\LieT\oplus\bigoplus_{\alpha\in R}{\LieG'}_{\alpha}$.
Applying $\ad H$ gives 
$ (\ad H) (B) = \sum_{\alpha\in R} \langle \alpha,H\rangle
B_{\alpha}$.
Combined with $B \in \Eig(\ad H; r)$, this
implies that
$$B=B_0\delta_{r0}+\sum_{\alpha\in R \atop \langle
  \alpha,H \rangle=r}
B_\alpha.$$
Now $(\Ad t) (Bz^{mr})$ is equal to
$$ \Big( B_0\delta_{r0} 
  +\sum_{\alpha\in R \atop \langle \alpha,H\rangle=r} 
  \alpha(t) B_\alpha \Big) z^{mr} 
= B_0\delta_{r0}z^{mr}+\sum_{\alpha\in R \atop \langle
  \alpha,H\rangle=r} B_\alpha =B.$$
\end{proof}

\begin{lemma}\label{L:torustrafo}
Let $m\in\DNplus$, $H \in \LieT$, and $\psi \in
\mathcal{X}^\cek(T)$ be such that
$\langle \alpha, \psi \rangle = m \langle \alpha, H \rangle$
for all $\alpha \in R_H^\DZ$.
Let $t=\psi (z\inv) \in T(K)$ 
and $Y = \widehat{z\delz}(t)t\inv \in \LieT(K)$.
Then we have $Y  \in \LieT$, $\langle \alpha, Y \rangle \in \DZ$ 
for all $\alpha \in R$, 
and $\langle \alpha, Y \rangle = -m \langle \alpha, H \rangle$ 
for all $\alpha \in R_H^\DZ$.
\end{lemma}
\begin{proof}
Let $d$ be the dimension of the torus $T$. We choose an isomorphism
from $T$ to $(G_\text{m})^d$ and identify \mbox{$T = (G_\text{m})^d$}.
For $i \in \{1,\dots, d\}$, we define the cocharacter $\delta_i:
G_\text{m} \ra T$ by $\delta_i(x)= (1,\dots, 1,x,1,\dots, 1)$.
Let the character $\epsilon_i:T \ra  G_\text{m}$ be defined by
$\epsilon_i(x_1,\dots, x_d) =  x_i$.
We have $\mathcal{X}^\cek(T)= \bigoplus_{i=1}^d \DZ \delta_i$
and $\mathcal{X}(T)= \bigoplus_{i=1}^d \DZ \epsilon_i$. 

Let $\psi =\sum f_i\delta_i$ for suitable $f_i \in \DZ$. 
From 
$ t = \psi(z\inv) 
= (z^{-f_1},\dots, z^{-f_d})$,
we deduce that
$Y = z \delz(t)t\inv = (-f_1,\dots, -f_d) \in \LieT$.
An arbitrary $\alpha \in R$ can be written as $\alpha=\sum a_i\epsilon_i$
for suitable $a_i \in \DZ$. Then we have
$\langle\alpha, Y\rangle = -\sum a_if_i \in \DZ$.
If $\alpha \in R_H^\DZ$, we obtain 
$m \langle \alpha, H \rangle = \langle \alpha, \psi \rangle = \sum
a_if_i = - \langle \alpha, Y \rangle$.
\end{proof}

\subsection{Transforming Regular Connections to Zero Standard
  Form}\label{SS:nullstandard}

\begin{definition}
Let $G$ be a linear algebraic group. We denote by 
$\LieGNull$ the set of all $X \in \LieG$ such that zero is the only rational
eigenvalue of $\ad X$.
The elements of $\LieGNull \dlogz$ are called {\bf connections in zero standard form}.
\end{definition}

\begin{theorem}\label{T:standardnull}
Let $G$ be a connected reductive linear algebraic group.
For every regular connection $A$, 
there exists a positive integer $n \in \DNplus$,
such that the pull-back  
connection $\incl{n}(A)$ is gauge equivalent
to a connection in zero standard form.
\end{theorem}
\begin{corollary}\label{c:regular-related-zero-standard}
  If $G$ is connected reductive, each regular connection is related to
  a connection in zero standard form. 
\end{corollary}
\begin{proof}
  Similar to the proof of Corollary~\ref{c:regular-related-standard}.
  Note that $\LieGNull$ is stable under multiplication by rational
  numbers.
\end{proof}

\begin{example}
   Let $G=\SL_2$. For $n \in
   \DNplus$, we define the regular $\SL_2$-connection
    \begin{equation*}
      B_n =  \tfrac1{2n}
      \begin{bmatrix}
          1 & 0 \\
          0 & -1
      \end{bmatrix}
      \dlogz.
    \end{equation*}
   We claim that $\incl{n}(B_n)$ is not gauge equivalent to
   a connection in zero standard form. Otherwise suppose that
   $$\incl{n}(B_n) 
   = 
   \begin{bmatrix}
     \tfrac12 & 0 \\
     0 & -\tfrac12
   \end{bmatrix}
   \dlogz$$
   is $\SL_2(K)$-equivalent to $X \dlogz$
   for some $X \in \LieSLzweiNull$. We may assume that 
   \mbox{$X \in \LieSLzweiNull$} has Jordan normal form.
   According to Proposition~\ref{p:slngleichglnaequi},
   there is $l \in \DZ$ such that 
   $$X = 
   \begin{bmatrix}
     \tfrac12 + l & 0 \\
     0 & -\tfrac12 -l
   \end{bmatrix}.$$
   Then $1+2l \in \DQ -\{0\}$ is a rational eigenvalue of $\ad X$. This
   contradicts our assumption $X \in \LieSLzweiNull$. 

   This example shows that there is no $n \in \DNplus$, such that
   for every regular $\SL_2$-connection $B$, the connection
   $\incl{n}(B)$ is gauge equivalent to a connection
   in zero standard form.
\end{example}

We prepare for the proof of Theorem~\ref{T:standardnull} by showing
some results for semisimple linear algebraic groups and
semisimple Lie algebras.

\begin{proposition}\label{P:jordancartan}
Let $G$ be a semisimple linear algebraic group, 
\mbox{$B \subset G$} a Borel subgroup,
and $T \subset B$ a maximal torus in $G$. 
The set of semisimple elements in $\LieB$ is equal to $\Ad(B)(\LieT)$.
\end{proposition}
\begin{proof}
Suppose that $H\in\LieB$ is a semisimple element. Let $\LieT'$ be a 
maximal toral subalgebra of $\LieG$, containing $H$. 
As $\LieG$ is semisimple, $\LieT'$ is a Cartan subalgebra of $\LieG$.
All Cartan subalgebras are conjugate under the adjoint group of $\LieG$. 
This adjoint group is the identity component of $\Ad(G)$. 
Therefore, we find $g \in G$ such that $(\Ad g) (\LieT)=\LieT'$.
Obviously, $\LieT'$ is the Lie algebra of $T'=gTg\inv$. 
The identity component $D$ of $T'\cap B$ is a torus.
We find 
$H \in \LieT' \cap \LieB = \Lie(T'\cap B) =\LieD$.
Let $S$ be a maximal torus in $B$ that contains $D$.
The maximal tori $T$ and $S$ in $B$ are conjugate under
$B$, so we find $b\in B$ such that $bTb\inv=S$. Now we see that
$(\Ad b)(\LieT)=\LieS\supset\LieD\ni H$.
\end{proof}
\begin{corollary}\label{K:jordancartan}
For every $X\subset \LieG$, there is a group element $g\in
G$ such that $(\Ad g) (X) \in \LieB$ and $(\Ad g) (X_{{\op{s}}}) \in
\LieT$.
\end{corollary}
\begin{proof}
Given $X \in \LieG$, 
let $\LieB'$ be a Borel subalgebra containing $X_{\op{s}}$ and $X_{\op{n}}$.
Since all Borel subalgebras are conjugate under the adjoint group, we
find an element $h \in G$ with
$(\Ad h)(\LieB')=\LieB$. By Proposition~\ref{P:jordancartan},
there is $b \in B$ such that
$(\Ad bh)(X_{\op{s}}) \in \LieT.$
We have $(\Ad bh)(X) \in \LieB$.
\end{proof}

\begin{lemma}\label{L:jordannilpotent}
Suppose that $\LieT$ is a Cartan subalgebra of a semisimple Lie
algebra $\LieG$, let $R=R(\LieG,\LieT)$ be the roots and choose
a system of positive roots $R^+ \subset R$.
Define $\LieU^+ = \bigoplus_{\alpha \in R^+}\LieG_\alpha$ and $\LieB =
\LieT \oplus \LieU^+$.
Let $X \in \LieB$ and $X=X_{\op{s}}+X_{\op{n}}$ be its Jordan decomposition
in $\LieG$. Then $X_{\op{s}} \in \LieT$ implies that $X_{\op{n}} \in \LieU^+$.
\end{lemma}
\begin{proof}
This follows from the root space decomposition, the nilpotency of $\ad
X_{\op{n}}$, and $\LieT^\ast= kR^+$.
\end{proof}

Once again, we use the notation introduced in Subsection~\ref{SS:standard}. 
In particular, $G$ is connected reductive and
$G'=(G,G)$ is semisimple.
Let $B'\subset G'$ be a Borel subgroup containing $T'$, and let
$R^+=R^+(B')\subset R$ be the corresponding positive roots. The set
${U'}^+ \subset B'$ of all unipotent elements of $B'$ is a closed
nilpotent connected subgroup.
For the corresponding Lie algebras, we have
${\LieU'}^+ = \bigoplus_{\alpha \in R^+}{\LieG'}_\alpha$
and
$\LieB' = {\LieT'} \oplus {\LieU'}^+$.

\begin{proposition}\label{P:jordananpassenred}
For every $X \in \LieG$, there is $g \in G$ such that 
$(\Ad g)(X_{\op{s}}) \in \LieT$ and
$(\Ad g)(X_{\op{n}}) \in {\LieU'}^+$.
\end{proposition}

\begin{proof}
This follows from Corollary~\ref{K:jordancartan},
Lemma~\ref{L:jordannilpotent} and the adaption to the reductive
situation.
\end{proof}

\begin{lemma}\label{L:eigenwertered}
Let $H \in \LieT$ and $N\in {\LieU'}^+$. The eigenvalues of
$\ad(H+N)$ are given by 
$\{\langle\alpha,H\rangle\mid \alpha \in R\cup\{0\}\}$.
\end{lemma}

\begin{proof}
The endomorphisms $\ad H$ and $\ad (H + N)$ have the same
characteristic polynomial.
\end{proof}

\begin{proof}[Proof of Theorem~\ref{T:standardnull}]
Let $A$ be a regular connection. 
According to Theorem~\ref{T:standard}, there are a positive integer $m
\in \DNplus$ and an element $X \in \LieG$,
such that $\incl{m}(A)$ is gauge equivalent to the connection
$C=X \dlogz$.

By Proposition~\ref{P:jordananpassenred}, we may assume that
$X_{\op{s}} \in \LieT$ and $X_{\op{n}} \in {\LieU'}^+$. 
Let
$$R_{X_{\op{s}}}^\DQ = \{\alpha \in R \mid \langle \alpha, X_{\op{s}}
\rangle \in \DQ\}$$
denote the roots rational on $X_{\op{s}}$. Choose
$l \in \DNplus$ such that
$\langle \alpha, lX_{\op{s}} \rangle \in \DZ$ for all 
$\alpha \in R_{X_{\op{s}}}^\DQ$.
Define $H=lX_{\op{s}}$ and $N=lX_{\op{n}}$.
If a root attains a rational value on $H$, this value is already
integral, which means $R_{H}^\DQ = R_{H}^\DZ$.

By Lemma~\ref{L:mHtuts}, there is $\phi \in
\mathcal{X}^\cek(T)$, such that 
$$\langle \alpha, \phi\rangle = m_H \langle\alpha, H\rangle
\quad\text{for all
$\alpha \in R_H^\DZ$}.$$ 
Define $t = \phi(z\inv) \in T(K)$.
Because
$$\alpha(t)=z^{-\langle\alpha,\phi\rangle}
=z^{-m_H\langle\alpha,H\rangle}
\quad\text{for all
$\alpha \in R_H^\DZ$},$$
and $m_HlX \in \Eig (\ad H; 0)$, Lemma~\ref{L:Adt} 
implies that $(\Ad t) (m_HlX)= m_HlX$.

Set $Y= \widehat{z \delz}(t)t\inv$. Lemma~\ref{L:torustrafo}
gives $Y \in \LieT$,
\begin{align}
\langle \alpha, Y \rangle \in& \;\DZ 
\quad\text{for all $\alpha \in R$, and}\label{Eq:Ywurzel}\\
\langle \alpha, Y \rangle =& -m_H \langle \alpha, H \rangle 
\quad\text{for all $\alpha \in R_H^\DZ$.}\label{Eq:Yganzewurzel}
\end{align}
We apply the gauge transformation $t$ to  
$\incl{(m_Hl)}(C) = m_Hl X\dlogz$ and get
$$ t\left[ \incl{(m_Hl)}(C) \right] 
= (\underbrace{m_HH + Y}_{\in \LieT} + \underbrace{m_HN}_{\in {\LieU'}^+})
   \dlogz.$$
Let $\alpha \in R$ with $\langle \alpha, m_HH + Y\rangle \in \DQ$.
We claim that $\langle\alpha, m_HH + Y\rangle = 0$.
Equation~\eqref{Eq:Ywurzel} shows that
$\alpha \in R_H^\DQ=R_H^\DZ$.
Then Equation~\eqref{Eq:Yganzewurzel} proves our claim.
Thus, by Lemma~\ref{L:eigenwertered}, $(m_HH+Y)+m_HN$ is an element of $\LieGNull$.
\end{proof}
\begin{corollary}\label{c:standardnull-semisimple-related}
Let $G$ be connected reductive with maximal torus $T \subset G$.
Then every connection of the form $X \dlogz$ with $X \in \LieT$ is related
to a connection in $\LieT \cap \LieGNull \dlogz$.
\end{corollary}
\begin{proof}
  If $X \in \LieT$, it is obvious from the proof of
  Theorem~\ref{T:standardnull} that there exists $n \in \DNplus$
  such that $\incl{n}(X \dlogz)$ is gauge equivalent
  to $X' \dlogz$ for some $X' \in \LieT \cap \LieGNull$.
  So $X \dlogz$ is related to $n\inv X' \dlogz$.
\end{proof}

\section{Regular \texorpdfstring{$\GL_n$}{GLn}-Connections}\label{S:dmoduln}

\subsection{Classification of Regular $\GL_n$-Connections}

For $a \in\DNplus$ and $x \in k$, we denote by
$\Einheit_a \in \End(k^a) = \Mat_{a}(k)$ 
the identity matrix and by 
$\Jordan(x, a) \in \Mat_{a}(k)$
the $a \times a$-Jordan block with diagonal entries equal to $x$.
For example, we have
$$ 
\Jordan(x,3)=
\begin{bmatrix}
x & 1 & 0 \\
0 & x & 1     \\
0 & 0 & x \\
\end{bmatrix}.
$$
Let $n \in \DN$. We denote the set of all $n\times n$-matrices in Jordan normal
form by ${\mathcal{J}_n}$.

\begin{definition} 
Let $X$, $Y \in {\mathcal{J}_n}$ be given,
\begin{align}
  \label{eq:XYjordan}
  \begin{split}
    X = & \blockdiag(\Jordan(x_1, a_1),\dots , \Jordan(x_r, a_r)),\\
    Y = & \blockdiag(\Jordan(y_1, b_1),\dots , \Jordan(y_s, b_s)).
  \end{split}
\end{align}
The matrices $X$ and $Y$ {\bf differ integrally}
(resp.\ {\bf rationally}) 
{\bf after block permutation}, if $s=r$ and there
is a permutation $\tau \in \Sym_r$ such that 
$$a_i = b_{\tau(i)} \quad \text{and}\quad 
x_i \equiv y_{\tau(i)} \mod \DZ
\quad\text{(resp. \negthickspace\negthickspace\negthickspace $\mod \DQ$)}\quad\text{
for all $i \in \{1, \dots, r\}$.}$$
\end{definition}

Now we can classify regular $\GL_n$-connections up to gauge
equivalence. Our
proof of the following theorem is more or less the same as that in \cite[\S3]{BV}.

\begin{theorem}[Classification of Regular $\GL_n$-Connections up to
  Gauge Equivalence]\label{T:GLnklass}
  The map
  $${\mathcal{J}_n} \sra \{\text{regular
    $\GL_n$-connections}\}/\GL_n(K), \quad X \mapsto [X \dlogz],$$
  is a surjection. For $X$, $Y \in {\mathcal{J}_n}$, the connections $X \dlogz$ and
  $Y \dlogz$ are $\GL_n(K)$-equivalent if and only if 
  $X$ and $Y$ differ integrally after block permutation.
\end{theorem}

\begin{proof}
  The surjectivity follows 
  from Theorem~\ref{T:standard} and the example $G=\GL_n$ in
  Examples~\ref{Ex:standard}. 
  Let $X$, $Y \in \mathcal{J}_n$ be given in the form
  \eqref{eq:XYjordan}.

  Assume that $X$ and $Y$ differ integrally after block
  permutation. Performing the block permutation by an element of
  $\GL_n(k)$ (or $\SL_n(k)$), we may assume that $Y$ has the form
  \begin{equation*}
    Y = \blockdiag(\Jordan(x_1+n_1, a_1),\dots , \Jordan(x_r+n_r, a_r))
  \end{equation*}
  for suitable $n_1,\dots, n_r \in \DZ$. We define
  \begin{equation}
    \label{eq:gform}
    g = \blockdiag\big(z^{n_1} \Einheit_{a_1}, \dots, z^{n_r}
    \Einheit_{a_r}\big) \in \GL_n(K)
  \end{equation}
  and deduce from Equation~\eqref{Eq:actionGLn} that $g\left[X
    \dlogz\right] = Y \dlogz$. 

  Assume now that $X \dlogz$ and $Y \dlogz$ are gauge equivalent. 
  We transform $X \dlogz$ by a gauge transformation $g$ of the
  form \eqref{eq:gform} for suitable $n_1, \dots, n_r \in \DZ$ and
  deal with $Y \dlogz$ similarly, and may so assume that
  \begin{equation}
    \label{eq:eigenvalue-cond}
    \lambda -\mu \in \DZ \Rightarrow \lambda =\mu \quad \text{for all
      $\lambda$, $\mu \in \{a_1,\dots a_r, b_1, \dots, b_s\}$.}
  \end{equation}
  Let $h \in \GL_n(K)$ with $h[Y \dlogz]=X \dlogz$. This
  implies $Xh-hY=z\delz(h)$. We write $h=\sum_{l \geq N} h_lz^l$ with
  $N \in \DZ$ and $h_l \in \Mat_n(k)$ and get $Xh_l-h_lY = lh_l$ for
  all $l\geq N$. Since the  eigenvalues of the linear map $\Mat_n(k)
  \ra \Mat_n(k)$, $A \mapsto XA-AY$ are given by $a_i-b_j$, for 
  $1\leq i \leq r$, $1 \leq j \leq s$, \eqref{eq:eigenvalue-cond}
  implies that $h
  =h_0 \in \GL_n(k)$. But then $h_0 Y h_0\inv =X$, so $X$ and
  $Y$ are conjugate and hence differ integrally (in fact, by $0$)
  after block permutation.
\end{proof}

\subsection{$D$-Modules and $\GL_n$-Connections}
We denote by $D_0= k[[z]]z\delz$ the subspace of derivations $\delta \in D$ with
$\delta(\calmax) \subset \calmax$.  

\begin{definition}[\cite{Manin}]
A {\bf $D$-module} is a $K$-vector space $M$ together with a map
$\alpha: D \times M \ra  M$, $(\delta, m) \mapsto \delta m =
\alpha(\delta, m)$,
that is $K$-linear in the first argument and additive in
the second one, and that satisfies  
$\delta(xm) = (\delta x) m + x (\delta m)$ for all $\delta \in D$, $x
\in K$, and $m \in M$. 
A map $\alpha$ as above is a
{\bf $D$-module structure} on $M$.
A {\bf morphism of $D$-modules}
$f:(M, \alpha) \ra (N, \beta)$ is a $K$-linear map $f: M
\ra N$ satisfying $f(\alpha(\delta, m)) = \beta (\delta, f(m))$ for all
$\delta \in D$, $m \in M$.

Let $M$ be a $D$-module. For $m \in M$, let $E(m)$ be the
smallest $\calO$-submodule of $M$ that contains $m$ and is
$D_0$-stable.
A $D$-module $M$ is {\bf Fuchsian}, 
if $E(m)$ is finitely generated as an $\calO$-module for all $m \in M$.
\end{definition}

Let $a \in \DNplus$ and $x \in k$. There is a unique $D$-module
structure $\alpha$ on $M = K^a=\bigoplus_{i=1}^{a}K e_i$ such that
$$\alpha(z\delz, e_i)
= x e_i + e_{i-1} \quad
\text{for all $i \in \{1,\dots a\}$},$$ 
where $e_0=0$.
We denote this $D$-module by $M^{x,a}$. It is easy to see that
$M^{x,a}$ is Fuchsian and indecomposable.

Fix $n \in \DN$. The group $\GL_n(K)$ acts on the set of all 
$D$-module structures on ${K^n}$ as follows: Given $g \in \GL_n(K)$
and a $D$-module structure $\alpha$, we define $g.\alpha$ by $(\delta,
w) \mapsto g(\alpha(\delta, g\inv(w)))$. Two $D$-module structures
$\alpha$ and $\beta$ are in the same orbit if and only if 
$({K^n}, \alpha)$ and $({K^n}, \beta)$ are isomorphic.

Let $A$ be a $\GL_n$-connection. 
If we evaluate
$A \in \CGLn = \Hom_{K}(D, \LieGL_n(K))$ 
at $\delta \in D$, we get an element
$A(\delta) \in \End_K({K^n}) = \Hom({k^n},{K^n})$.
Then
$$\alpha_A (\delta, x\otimes v) = \delta(x)\otimes v - x
A(\delta)(v),$$ 
where $x \in K$, $v \in k^n$, defines a $D$-module structure
$\alpha_A$ on ${K^n}=K\otimes {k^n}$. 

We omit the easy proof of
\begin{proposition}\label{P:dmodulGLnzshg}
The map 
$$
  \{\text{$\GL_n$-connections}\} \sira 
  \{\text{$D$-module structures on ${K^n}$}\}, \quad
  A  \mapsto  \alpha_A,$$
is a $\GL_n(K)$-equivariant bijection and induces a bijection between
regular connections and Fuchsian $D$-module structures.
The regular $\GL_n$-connection
$$\blockdiag(-\Jordan(x_1, a_1),\dots, -\Jordan(x_r, a_r)) \dlogz$$
corresponds to the Fuchsian $D$-module
$M^{x_1,a_1}\oplus\dots \oplus M^{x_r,a_r}$.
\end{proposition}

We use Proposition~\ref{P:dmodulGLnzshg} in order to translate
Theorem~\ref{T:GLnklass} in the language of $D$-modules and obtain

\begin{theorem}[cf.~{\cite[Theorem~4]{Manin}}]\label{T:fuchsklass}
Every finite-dimensional Fuchsian $D$-module is a direct
sum of indecomposable Fuchsian $D$-modules.
The summands are unique up to permutation and isomorphism.

The map $(x,a) \mapsto M^{x,a}$ induces a bijection from $k/\DZ
\times \DNplus$ to the set of isomorphism classes of
finite-dimensional indecomposable Fuchsian $D$-modules. 
\end{theorem}

\section{Regular
  \texorpdfstring{$\SL_n$}{SLn}-Connections}\label{S:klassifikation}
Let $\mathcal{J}(\LieSL_n)=\LieSL_n \cap \mathcal{J}_n$ and
$\mathcal{J}(\LieSLnNull)=\LieSLnNull \cap \mathcal{J}_n$.
This section is organized as follows.
First we describe the set $\Rel(X\dlogz)/\SL_n(K)$ of relatives up to
gauge equivalence, for $X \in
\mathcal{J}(\LieSLnNull)$ (Bijection~\eqref{eq:kombverwandt}). We deduce
from this
that the set of regular $\SL_n$-connections is 
$\bigcup\Rel(X \dlogz)$, where $X$ ranges over $\mathcal{J}(\LieSLnNull)$
(Proposition~\ref{p:alleregulaer}, Corollary~\ref{K:alleverwandtenregulaer}).
Then we establish the classification up to relationship
(Theorem~\ref{t:sln-class-rel})
and up to gauge equivalence (Theorem~\ref{S:slnklassifik}). 
We conclude with a slightly different view on this classification
(Remark~\ref{rem:nice-class}) and
explain the example $\SL_2$ (Example~\ref{ex:SL2}).

Let $X \in \LieSLnNull$. Recall from Proposition~\ref{P:verwandteh1}
the bijection  
\begin{equation}
  \label{eq:rel-H1}
  \Rel(X\dlogz)/\SL_n(K) \sira \Ho^1(K; X\dlogz).
\end{equation}
Since $\Ho^1(K; X\dlogz)$ is the direct
limit of the $\Ho^1(\Gamma_l; \SL_n(L)_{\incl{l}(X \dlogz)})$, for $l \in
\DNplus$ and $\incl{l}:K \hra K=L$, we are interested in the
stabilizer of $\incl{l}(X \dlogz)=lX\dlogz$ in $\SL_n(L)$.
As $\LieSLnNull$ is stable under multiplication by rational
numbers, Proposition~\ref{P:stabisln} shows that
\begin{equation*}
  \SL_n(L)_{lX \dlogz}=\Z_{\SL_n}(lX)=\Z_{\SL_n}(X) \subset \SL_n(k). 
\end{equation*}
In particular, the action of the Galois group $\Gamma_l$ on
$\SL_n(L)_{lX \dlogz}$ is trivial.
\begin{proposition}\label{P:stabisln}
Let $X \in \LieSLnNull$. Then $\SL_n(K)_{X \dlogz}=\Z_{\SL_n}(X)$,
where $\Z_{\SL_n}(X)$ is the centralizer of $X$ 
in ${\SL_n}={\SL_n}(k)$ under the adjoint action.
\end{proposition}
\begin{proof}
  The inclusion $\Z_{\SL_n}(X) \subset {\SL_n}(K)_{X \dlogz}$ is obvious.
  Let \mbox{$g \in \SL_n(K)_{X \dlogz}$}. As
  $\SL_n(K) \subset \GL_n(K) \subset \Mat_{n}(K)$, we find 
  $N \in \DZ$ and $g_i \in \Mat_{n}(k)$ such that
  $g = \sum_{i\ge N} g_i z^i$.
  From $g[X \dlogz]= X \dlogz$ and Equation
  \eqref{Eq:actionGLn} we get $ Xg-gX = z\delz(g)$, or,
  equivalently, 
  $(\ad_{\LieGL_n} X) (g_i) = ig_i$ for all $i \ge N$.
  But then $X \in \LieSLnNull \subset \LieGLnNull$
  implies that $g_i= 0$ for $i \neq 0$, in other words,
  $g = g_0 \in \Z_{\SL_n}(X)$.
\end{proof}
Recall that $\ol{K}=\bigcup_{m\in\DNplus}k((z^{1/m}))$
is an algebraic closure of $K=k((z))$.
For $l \in \DNplus$, we view the field extension $\incl{l}: K \hra
K=L$ as a subextension of $K \subset \ol{K}$ via the embedding
$K=L \hra \ol{K}$, $f(z) \mapsto f(z^{1/l})$.
The Galois group
$\Gal(\ol{K}/K)$ is isomorphic to the procyclic group
$\widehat{\DZ}$.
For the rest of this section, we fix
a procyclic generator $\gamma$ of $\Gal(\ol{K}/K)$.
For $l \in \DNplus$, 
the Galois group $\Gamma_l=\Gal(\incl{l})$ is generated by $\gamma|_{L}$. 

Let $X \in \LieSLnNull$ and $l \in \DNplus$. 
Since $\Gamma_l$ acts trivially on $\SL_n(L)_{lX \dlogz}$,
the map 
\begin{equation}\label{z1isotel}
\Z^1(\Gamma_l; {\SL_n}(L)_{lX \dlogz}) \xrightarrow[\gamma]{\sim}
\{\text{$l$-torsion-elements in $\Z_{\SL_n}(X)$}\}, \quad
p \mapsto p_{\gamma|_L},
\end{equation}
is bijective. It depends on $\gamma$. We
indicate this here and in the following by putting a $\gamma$ at the
corresponding arrow. 
This Bijection \eqref{z1isotel} induces a bijection
\begin{equation*}%\label{h1isotel}
\Ho^1(\Gamma_l; {\SL_n}(L)_{lX \dlogz}) \xrightarrow[\gamma]{\sim} 
\{\text{conj.\ classes of $l$-torsion-elts.\ in $\Z_{\SL_n}(X)$}\}.
\end{equation*}
We use Proposition~\ref{P:H1fieldinclusion} and pass to the direct
limit. We obtain for $X \in \LieSLnNull$ a bijection 
\begin{equation}
  \label{eq:h1kktorsion}
 {\Ho^1(K; X \dlogz)} \xrightarrow[\gamma]{\sim} 
\{\text{conj.\ classes of torsion elements in $\Z_{\SL_n}(X)$}\}.
\end{equation}

Suppose that $X \in \mathcal{J}(\LieSL_n)$.
Let $T_X \subset \Z_{\GL_n}(X)$ be the diagonal maximal torus and
$W_X$ be the Weyl group as in Theorem~\ref{T:zentralisator}. 
The Weyl group $W_X$ stabilizes $D_X = T_X \cap \SL_n$.

\begin{proposition}\label{P:TEDX}
  For $X \in \mathcal{J}(\LieSL_n)$, the inclusion $D_X \hra
  Z_{\SL_n}(X)$ induces a bijection
  \begin{equation*}
    \{\text{torsion elements in $D_X$}\}/W_X \sira 
    \{\text{conj.\ classes of torsion elts.\ in $\Z_{\SL_n}(X)$}\}.
  \end{equation*}
\end{proposition}
\begin{proof}
  This follows from 
  Theorem~\ref{T:zentralisator}~\eqref{En:heKk} and the fact that,
  in characteristic zero, every element of finite order is semisimple.
\end{proof}

Assume now that $X \in \mathcal{J}(\LieSLnNull)$.
We combine 
Bijections~\eqref{eq:rel-H1}, \eqref{eq:h1kktorsion} and Proposition~\ref{P:TEDX}
in order to get the bijection
\begin{equation}
  \label{eq:kombverwandt}
  {\{\text{torsion elements in $D_X$}\}/W_X}  \xrightarrow[\gamma]{\sim}
  \Rel(X \dlogz)/\SL_n(K).  
\end{equation}
We denote this map by $\delta \mapsto \left[(X \dlogz)^\delta\right]$
and describe it explicitely in the proof of 
\begin{proposition}\label{p:alleregulaer}
  All relatives of $X \dlogz$ are regular, for
  $X \in \mathcal{J}(\LieSLnNull)$. 
\end{proposition}
\begin{proof}
Write $X \in \mathcal{J}(\LieSLnNull)$ as
$$X=\blockdiag(\Jordan(x_1, a_1),\dots, \Jordan(x_r, a_r))$$
for suitable $x_i \in k$ and $a_i \in \DNplus$.
Given a torsion element $d \in D_X$ we explain now how to construct 
an $\SL_n$-connection in the orbit
$\left[(X \dlogz)^{W_X d}\right]$. As this connection will be regular,
this proves the proposition.
For $l \in \DNplus$, we view $\incl{l}$ as the field extension
$K=k((z)) \hra k((z^{1/l}))$. Let $\omega_l$ be the primitive $l$-th root of
unity such that
$\gamma(z^{1/l})= \omega_l z^{1/l}$.

Let $d \in D_X$ be a torsion element. We find $l \in \DNplus$ 
and $j_1,\dots, j_r \in \DN$ such that
$$d= \blockdiag\big(\omega_l^{j_1}\Einheit_{a_1}, \dots, \omega_l^{j_r}\Einheit_{a_r}\big).$$ 
Let $\omega=\omega_l$, $\zeta=z^{1/l}$, and $\Sigma=\sum_{s=1}^r
{j_sa_s} \in \DN$.  
As $1=\det (d)=\omega^\Sigma$, we see that 
$\Sigma$ is divisible by $l$. Define
$$g=\blockdiag\big(\zeta^{j_1}\Einheit_{a_1}, \dots,
\zeta^{j_{r-1}}\Einheit_{a_{r-1}}, 
  \zeta^{j_r}, \dots, \zeta^{j_r}, \zeta^{j_r-\Sigma}\big) \in
  \SL_n((\zeta)).$$ 
Now $\omega^\Sigma=1$ implies that $d=g\inv \gamma(g)$.
Therefore, for any $m \in \DNplus$, we get
$$  d^m =  d \gamma(d) \gamma^2(d)\cdots \gamma^{m-1}(d)
        =  g\inv \gamma^m(g).$$ 
This means that $d$, regarded as an element of $\Z^1(\Gamma_l;
\SL_n((\zeta))_{lX \dlogzeta})$ via Bijection~\eqref{z1isotel}, 
is cohomologous to the trivial 1-cocycle
in $\Z^1(\Gamma_l; \SL_n((\zeta)))$.
From the proof of Theorem~\ref{T:H1Kform} results that the connection
$g\left[lX\dlogzeta\right]$
is invariant under the Galois group $\Gamma_l$.
We define
$$(X \dlogz)^d = \blockdiag\Big(\Jordan(x_1+\tfrac{j_1}{l},
    a_1), \dots, \Jordan(x_{r-1}+\tfrac{j_{r-1}}{l},
    a_{r-1}), C\Big) \dlogz,$$
where $C \in \Mat_{a_r}(K)$ is given by
$$ C= 
\begin{bmatrix}
x_r+\tfrac{j_r}{l}   & 1     & 0      & \dots      &       &    0   \\
0      & x_r+\tfrac{j_r}{l}   & \ddots & \ddots      &       &       \\
\vdots  & \ddots      & \ddots &  1     &   0    &   \vdots    \\
      &       & 0      & x_r+\tfrac{j_r}{l}   & 1     &   0    \\
      &       &       & 0      & x_r+\tfrac{j_r}{l}   & z^{\frac{\Sigma}{l}} \\
0      &       &       & \dots      & 0      &
x_r+\tfrac{j_r-\Sigma}{l}   
\end{bmatrix}
$$
If $a_r=1$, this is to be interpreted as $C=x_r+\tfrac{j_r-\Sigma}{l}$.   
As $\tfrac{\Sigma}{l} \in \DN$ is a nonnegative integer, $(X
\dlogz)^d$ is
regular. It is easy to verify that
$\incl{l}\left((X \dlogz)^d\right)=g\left[lX \dlogzeta\right]$.
We conclude that $(X \dlogz)^d$ is a connection in 
$\left[(X \dlogz)^{W_X d}\right]$. 
\end{proof}

\begin{corollary}\label{c:diagonal-rel-equi}
  If $X\in \LieSL_n$ is a diagonal matrix, each connection related to
  $X \dlogz$ is gauge equivalent to a connection of the form $Y
  \dlogz$ with $Y \in \LieSL_n$ a diagonal matrix.
\end{corollary}
\begin{proof}
  By Corollary~\ref{c:standardnull-semisimple-related} we may assume
  that $X \in \LieSLnNull$ is diagonal. Then our claim follows from the description
  of Bijection~\eqref{eq:kombverwandt} in the above proof.
\end{proof}

\begin{corollary}\label{K:alleverwandtenregulaer}
Every regular $\SL_n$-connection is related to $X \dlogz$, for some $X
\in \mathcal{J}(\LieSLnNull)$.
All relatives of a regular $\SL_n$-connection are regular.
An $\SL_n$-connection $A$ is regular if and only if there is 
$l \in \DNplus$ such that the connection $\incl{l}(A)$ is regular.
\end{corollary}
\begin{question}\label{q:grelatregulaer}
Are all relatives of a regular $G$-connection regular, if
$G$ is an arbitrary linear algebraic group?
\end{question}
\begin{proof}%[Proof of Corollary~\ref{K:alleverwandtenregulaer}]
The first claim follows from
Corollary~\ref{c:regular-related-zero-standard}, and then
the second claim is a consequence of Proposition~\ref{p:alleregulaer}.
If $\incl{l}(A)$ is regular, we have just seen that it is related to
$X \dlogz$ with $X \in \LieSL_n$. 
So $A$ is related to the regular connection $l\inv X\dlogz$ and
therefore regular.
\end{proof}
\begin{remark}\label{r:corollarygln}
Very similar arguments prove that Corollary
\ref{K:alleverwandtenregulaer} with $\SL_n$ replaced by
$\GL_n$ is true.
\end{remark}

\begin{proposition}\label{p:slngleichglnaequi}
For $X$, $Y \in \mathcal{J}(\LieSL_n)$,
the following are equivalent:
\begin{enumerate}
\item\label{SLequi} $X \dlogz$ and $Y \dlogz$ are $\SL_n(K)$-equivalent.
\item\label{GLequi} $X \dlogz$ and $Y \dlogz$ are $\GL_n(K)$-equivalent.
\item\label{differZ} $X$ and $Y$ differ integrally after block permutation.
\end{enumerate}
\end{proposition}

\begin{proof}
The implication \eqref{SLequi} $\Rightarrow$ \eqref{GLequi} is
obvious, and
\eqref{GLequi} $\Rightarrow$ \eqref{differZ} follows
from Theorem~\ref{T:GLnklass}.
In the proof of Theorem~\ref{T:GLnklass} we proved the implication
\eqref{differZ} $\Rightarrow$ \eqref{GLequi}. But we actually 
showed \eqref{differZ} $\Rightarrow$ \eqref{SLequi}: If the traces
of $X$ and $Y$ vanish, the element $g$ defined in
Equation~\eqref{eq:gform} is an element of $\SL_n(K)$.
\end{proof}

\begin{theorem}[Classification of Regular $\SL_n$-Connections up to
  Relationship]\label{t:sln-class-rel} 
  The map $X \mapsto X\dlogz$ induces a surjection
  \begin{equation*}
    \mathcal{J}(\LieSL_n) \sra \left\{\text{regular
        $\SL_n$-connections}\right\}/\text{relationship}.
  \end{equation*}
  For $X$, $Y \in \mathcal{J}(\LieSL_n)$,
  the connections $X \dlogz$ and $Y \dlogz$ are related if and only if $X$ and $Y$
  differ rationally after block permutation.
\end{theorem}

\begin{proof}
  Our map is surjective by Corollary~\ref{K:alleverwandtenregulaer}.
  The second statement follows from Proposition~\ref{p:slngleichglnaequi} and the fact that
  the matrices $l\Jordan(x,a)$ and $\Jordan(lx,a)$ are
  $\SL_a(k)$-conjugate, for $l \in \DNplus$.
\end{proof}

Let $X \dlogz$ and $Y \dlogz$ be two related
$\SL_n$-connections with $X$, $Y \in \mathcal{J}(\LieSLnNull)$. 
From Bijection~\eqref{eq:kombverwandt}, we conclude that there is a unique map
$\can_{YX}$ such that the diagram
\begin{equation}
  \label{eq:canXY-diagram}
  \xymatrix@R-10pt{
    {\left\{\text{torsion elements in $D_X$} \right\}/W_X
      \ar@{->}[r]^-{\sim}_-{\gamma}} \ar[d]_-{\can_{YX}}^-\sim &  
    % {\left\{\text{relatives of $X \dlogz$} \right\}}
    \Rel(X \dlogz)/\SL_n(K) \ar@{=}[d]\\
    {\left\{\text{torsion elements in $D_Y$} \right\}/W_Y 
      \ar@{->}[r]^-\sim_-{\gamma}} &  
    % {\left\{\text{relatives of $Y \dlogz$} \right\}}
    \Rel(Y \dlogz)/\SL_n(K)
  }
\end{equation}
commutes. This map $\can_{YX}$ can be described explicitly, see
\cite{OSdiplom}.

\begin{theorem}[Classification of Regular $\SL_n$-Connections up to
  Gauge Equivalence]\label{S:slnklassifik} 
Let $\gamma$ be a procyclic generator of $\Gal(\ol{K}/K)$ and 
\begin{equation*}
  \coprod_{X \in \mathcal{J}(\LieSLnNull)}
  \{\text{torsion-elts.\ in $D_X$}\}/W_X \underset{\gamma}{\twoheadrightarrow}
  \left\{\text{regular $\SL_n$-connections}\right\}/\SL_n(K)
\end{equation*}
be the map induced by the maps \eqref{eq:kombverwandt},
$\delta \mapsto \left[(X \dlogz)^\delta\right]$.
Then this map is surjective, and we have 
$\left[(X \dlogz)^\delta\right] =\left[(Y \dlogz)^\epsilon\right]$ 
if and only if 
$X$ and $Y$ differ rationally after block permutation and
$\can_{YX}(\delta)= \epsilon$. 
\end{theorem}

\begin{remark}\label{rem:slnklassifik} 
The sets $\{\text{torsion elements in $D_X$}\}/W_X$ are
easy to
describe. ``Differing rationally after block permutation'' is an equivalence
relation on $\mathcal{J}(\LieSLnNull)$. 
By choosing a complete system of representatives for this
relation, and by using the explicit description of the map 
$\delta \mapsto \left[(X \dlogz)^\delta\right]$ given in the proof of
Proposition~\ref{p:alleregulaer}, 
Theorem~\ref{S:slnklassifik} enables us 
to give a list of all regular
$\SL_n$-connections up to $\SL_n(K)$-equivalence.
\end{remark}
\begin{proof}
Proposition~\ref{p:alleregulaer}, Bijection~\eqref{eq:kombverwandt} and
Corollary~\ref{K:alleverwandtenregulaer} 
show that our map is well defined and surjective.
The remaining claim follows from
Theorem~\ref{t:sln-class-rel} 
and Diagram~\eqref{eq:canXY-diagram}.
\end{proof}

\begin{remark}\label{rem:nice-class}
We now explain a nice partial classification of regular connections up
to gauge equivalence.
We associate to the standard diagonal Cartan subalgebra 
$\LieT \subset \LieSL_n$ the coroots $R^\cek$ and the Weyl group
$W$. This Weyl group acts naturally on $\LieT$ and stabilizes the 
subgroups $\DZ R^\cek$ and $\DQ R^\cek$. The groups
$W^\DZ = \DZ R^\cek \semidirect W$ and
$W^\DQ = \DQ R^\cek \semidirect W$ act on $\LieT$ by $(a,w).H=a+wH$.
Two elements
of $\LieT$ are in the same $W^\DZ$-orbit
(resp.\ $W^\DQ$-orbit) if and only if they differ integrally
(resp.\ rationally) after block permutation.
Let $\mathcal{N}_n= \mathcal{J}(\LieSL_n)-\LieT$. Consider the commutative
diagram 
\begin{equation}
  \label{eq:LieTmodW}
  \xymatrix@R-10pt{
    \LieT/W^\DZ \ar@{^{(}->}[r] \ar@{-{>>}}[d] &
    \left\{\text{regular
        $\SL_n$-connections}\right\}/\SL_n(K) \ar@{-{>>}}[d]^\pi \\
    % \mathcal{N}_n \ar[l] \ar@{=}[d] \\
    \LieT/W^\DQ \ar@{^{(}->}[r] &
    \left\{\text{regular
        $\SL_n$-connections}\right\}/\text{relationship} &
    \mathcal{N}_n \ar[l]_-{\nu} \\
  }
\end{equation}
with obvious vertical maps.
The horizontal maps are induced by $X \mapsto X \dlogz$. The
horizontal maps on the left are well-defined and injective by  
Proposition~\ref{p:slngleichglnaequi} and
Theorem~\ref{t:sln-class-rel}. They yield a partial classification.
In the lower row of Diagram~\eqref{eq:LieTmodW} the images of the
horizontal maps are complementary by Theorem~\ref{t:sln-class-rel}.
From Corollary~\ref{c:diagonal-rel-equi} follows that the image of
the upper horizontal map is the complement of
$\pi\inv(\nu(\mathcal{N}_n))$.
\end{remark}

\begin{example}\label{ex:SL2}
  We restrict now to the case $n=2$. Then $\mathcal{N}_2 =
  \{\Jordan(0,2)\}$, and $\nu(\mathcal{N}_2)$ consists of one element,
  namely $\Rel(\Jordan(0,2)\dlogz)/\text{relationship}$. 
  Its inverse image under $\pi$ is 
  $\Rel(\Jordan(0,2)\dlogz)/\SL_n(K)$. 
  From the description of the map \eqref{eq:kombverwandt} 
  in the proof of Proposition~\ref{p:alleregulaer} we see that this set 
  has precisely two elements, namely the orbits of the
  two connections 
  (cf.~example $G=\SL_2$ in Examples~\ref{Ex:standard})   
  \begin{equation*}
    \begin{bmatrix}
      0 & 1\\
      & 0
    \end{bmatrix}
    \dlogz
    \quad
    \text{and}
    \quad
    \begin{bmatrix}
      \tfrac 12 & z\\
      & -\tfrac 12
    \end{bmatrix}
    \dlogz.
  \end{equation*}
\end{example}

\section{Fuchsian Connections}\label{S:fuchszshg}

Let $G$ be a linear algebraic group and $\rho: G \ra
\GL(V)$ be a (rational)
representation of $G$ in a finite-dimensional vector space $V$. 
If $A$ is a $G$-connection, 
$\rho(A)=\rho_{\ast}(A) \in \CGL(V)$ is a
$\GL(V)$-connection and corresponds 
to a $D$-module structure $\alpha_{\rho(A)}$ on $K\otimes V$
(cf.~Proposition~\ref{P:dmodulGLnzshg}).

\begin{definition}
A connection $A$ is {\bf Fuchsian}
if for every finite-dimensional representation $\rho: G \ra
\GL(V)$ the $D$-module $(K\otimes V, \alpha_{\rho(A)})$ is Fuchsian.  
\end{definition}

\begin{remark}\label{B:fuchsaequidef}
  According to Proposition~\ref{P:dmodulGLnzshg}, a connection
  $A$ is Fuchsian if and only if for every finite-dimensional
  representation $\rho: G \ra \GL(V)$ the connection
  $\rho(A)$ is regular.
\end{remark}

\begin{proposition}
  Let $G$ be a linear algebraic group.
  Every regular connection is Fuchsian.
  For $G=\GL_n$ or \mbox{$G=\SL_n$}, every Fuchsian
  connection is regular.
\end{proposition}
\begin{question}\label{q:Fuchsgleichreg}
  Do the notions of Fuchsian and regular connection
  coincide for every linear algebraic group?
\end{question}
\begin{remark}
  Using Remarks~\ref{r:corollarygln} and \ref{B:fuchsaequidef}, it is
  easy to see that all relatives of a Fuchsian connection are Fuchsian.
  This shows that the answer ``yes'' to
  Question~\ref{q:Fuchsgleichreg} implies the same answer to Question~\ref{q:grelatregulaer} 
\end{remark}

\begin{proof}
The first claim and the second one for $G=\GL_n$ are obvious.
Let $A$ be a Fuchsian $\SL_n$-connection. Let \mbox{$\rho:\SL_n
\hra \GL_n$} be the standard representation of $\SL_n$.
There are $g \in \GL_n(K)$ and $X(z) \in
\LieGL_n[[z]]$ such that $g[\rho(A)]= X(z) \dlogz$.
Consider the field extension $\incl{n}: K \hra K=N$.
Let $f \in N$ be an $n$-th root of
$\incl{n}(\det(g\inv))$.
Then $h = f \incl{n}(g)$ is an element of $\SL_n(N)$,
and we have
$$ \rho(h[\incl{n}(A)]) =  f\Eins\left[ \incl{n}(X(z)\dlogz) \right]
                =  \left( nX(z^n) + z\delz(f) f\inv
                  \Eins\right) \dlogz.$$
It is obvious that $z\delz(f) f\inv \in k[[z]]$.
But then $h[\incl{n}(A)]$ is regular, and 
Corollary~\ref{K:alleverwandtenregulaer}
shows that $A$ is regular.
\end{proof}

\begin{appendix}
\section{Semisimple Conjugacy Classes}\label{App:cent}
Let $n \in \DN$ and $X \in \LieGL_n = \End({k^n})$.
For $\lambda \in k$ and $i \in \DN$, define
$$E_\lambda^i = \Ker (X - \lambda) \cap \Bild (X-\lambda)^i.$$
Every $g \in \Z_{\GL_n}(X)$ stabilizes all $E_\lambda^i$.
Thus $g$ induces maps
$g|_{E_\lambda^i} \in \GL(E_\lambda^i)$ and
$\ol{g|_{E_\lambda^i}} \in \GL(E_\lambda^i / E_\lambda^{i+1})$.
The following theorem gives an explicit description
of the semisimple conjugacy classes in $\Z_{\GL_n}(X)$.
\begin{theorem}\label{T:zentralisator}
  Let $n \in \DN$ and $X \in \LieGL_n = \Mat_n(k)$ be in
  Jordan normal form. Let $T_n \subset \GL_n$ be the standard diagonal
  torus. Then we have the following:
  \begin{enumerate}
  \item\label{En:torus} $T = T_X = T_n \cap \Z_{\GL_n}(X)$
    is a maximal torus in $\Z_{\GL_n}(X)$.
  \item\label{En:torusiso} The homomorphism 
    $$ \pi: \Z_{\GL_n}(X) \sra   
    \prod_{\lambda \in k \atop i \in \DN}\GL(E_\lambda^i /
    E_\lambda^{i+1}), \quad
    g \mapsto   \left(\ol{g|_{E_\lambda^i}}\right)_{\lambda, i},$$
    is surjective. It induces, by restriction, an isomorphism
    $ \pi|_{T}: T \sira  \pi(T)$,
    and $\pi(T)$ is a maximal torus in $ \prod \GL(E_\lambda^i / E_\lambda^{i+1})$.
  \item\label{En:heKk} The Weyl group $W_X$ associated 
    to the torus $\pi(T)$ in 
    $ \prod \GL(E_\lambda^i / E_\lambda^{i+1})$
    acts via $\pi|_T$ on T, and the inclusion $T_X \hra \Z_{\GL_n}(X)$
    induces a bijection
    \begin{equation*}
      T_X/W_X \sira \{\text{semisimple conj.\ classes in
        $\Z_{\GL_n}(X)$}\}.
    \end{equation*}
  \end{enumerate}
\end{theorem}
The proof is left to the reader. It can be found in \cite{OSdiplom}.
\end{appendix}

\bibliographystyle{amsalpha}

\begin{thebibliography}{Man65}

\bibitem[BV83]{BV}
Donald~G. Babbitt and Veeravalli~S. Varadarajan, \emph{Formal reduction theory
  of meromorphic differential equations: a group theoretic view}, Pacific J.
  Math. \textbf{109} (1983), no.~1, 1--80. \MR{86b:34010}

\bibitem[DG70]{DG}
Michel Demazure and Pierre Gabriel, \emph{Groupes alg{\'e}briques. {T}ome {I}:
  {G}{\'e}om{\'e}trie alg{\'e}brique, g{\'e}n{\'e}ralit{\'e}s, groupes commutatifs}, Masson
  \& Cie, {\'E}diteur, Paris, 1970. \MR{46 \#1800}

\bibitem[Lan52]{Lang}
Serge Lang, \emph{On quasi algebraic closure}, Ann. of Math. (2) \textbf{55}
  (1952), 373--390. \MR{13,726d}

\bibitem[Man65]{Manin}
Juri~I. Manin, \emph{Moduli fuchsiani}, Ann. Scuola Norm. Sup. Pisa (3)
  \textbf{19} (1965), 113--126. \MR{31 \#4815}

\bibitem[Sch03]{OSdiplom}
Olaf~M. Schn{\"u}rer, \emph{{R}egul{\"a}re {Z}usammenh{\"a}nge in trivialen
  algebraischen {$G$}-{H}auptfaserb{\"u}ndeln {\"u}ber der infinitesimalen
  punktierten {K}reisscheibe}, Diplomarbeit, Freiburg (2003),
  {\href{http://www.freidok.uni-freiburg.de/volltexte/1477/}{http://www.freidok.uni-freiburg.de/volltexte/1477/}}.

\bibitem[Ser68]{Serre}
Jean-Pierre Serre, \emph{Corps locaux}, Hermann, Paris, 1968, Deuxi{\`e}me
  {\'e}dition, Publications de l'Universit{\'e} de Nancago, No. VIII. \MR{50
  \#7096}

\bibitem[Ser97]{Serregaloiscoho}
\bysame, \emph{Galois cohomology}, Springer-Verlag, Berlin, 1997, Translated
  from the French by Patrick Ion and revised by the author. \MR{98g:12007}

\bibitem[Spr79]{SpringerAMS}
Tonny~A. Springer, \emph{Reductive groups}, Automorphic forms, representations
  and $L$-functions (Proc. Sympos. Pure Math., Oregon State Univ., Corvallis,
  Ore., 1977), Part 1, Proc. Sympos. Pure Math., XXXIII, Amer. Math. Soc.,
  Providence, R.I., 1979, pp.~3--27. \MR{80h:20062}

\bibitem[Spr98]{Springerneu}
\bysame, \emph{Linear algebraic groups}, second ed., Progress in Mathematics,
  vol.~9, Birkh{\"a}user Boston Inc., Boston, MA, 1998. \MR{99h:20075}

\bibitem[Ste65]{Steinberg}
Robert Steinberg, \emph{Regular elements of semisimple algebraic groups}, Inst.
  Hautes {\'E}tudes Sci. Publ. Math. (1965), no.~25, 49--80. \MR{31 \#4788}

\end{thebibliography}

\providecommand{\bysame}{\leavevmode\hbox to3em{\hrulefill}\thinspace}
\providecommand{\MR}{\relax\ifhmode\unskip\space\fi MR }
\providecommand{\MRhref}[2]{%
  \href{http://www.ams.org/mathscinet-getitem?mr=#1}{#2}
}
\providecommand{\href}[2]{#2}

\end{document}